\documentclass[runningheads]{llncs}
\usepackage{graphicx}
\usepackage[]{algorithm2e}
\SetKwInOut{Parameter}{Parameters}
\usepackage{amsmath,amssymb,amsfonts}
\usepackage{graphicx}
\usepackage{textcomp}
\usepackage{subcaption}
\usepackage{hyperref}
\begin{document}
\title{Automated differential equation solver based on the parametric approximation optimization}
%
%
\author{Alexander Hvatov \and Tatiana Tikhonova}
\authorrunning{A. Hvatov \and T. Tikhonova}
\titlerunning{Automated differential equation solver}
%
\institute{ITMO University, Kronsersky pr. 49, 197101, St. Petersburg, Russia
\email{alex\_hvatov@itmo.ru}
}
\maketitle              
\begin{abstract}
The numerical methods for differential equation solution allow obtaining a discrete field that converges towards the solution if the method is applied to the correct problem. Nevertheless, the numerical methods have the restricted class of the equations, on which the convergence with a given parameter set or range is proved. Only a few "cheap and dirty" numerical methods converge on a wide class of equations without parameter tuning with the lower approximation order price. The article presents a method that uses an optimization algorithm to obtain a solution using the parameterized approximation. The result may not be as precise as an expert one. However, it allows solving the wide class of equations in an automated manner without the algorithm's parameters change.

\end{abstract}
\section{Introduction}
\label{sec:intro}

Differential equations: ordinary (ODE) and partial differential equations (PDE) are classical ways to express the physical laws. The amount of the differential equations analyzed in the mathematical physics domain is limited, at least by the number of variational principles in physics. Modern data-driven methods \cite{maslyaev2021partial,rudy2017data,long2019pde}  may be considered as a source of yet not analyzed equations. During the equation discovery process, we obtain models with an arbitrary form, which may be a challenge for modern equation solvers.

An expert can always find a numerical solution for the obtained equation. However, a precise numerical solution of an equation of an unknown type is a challenging task. There are many expert systems and solvers for ODE systems that are alleviating the expert's work \cite{rackauckas2019confederated}. However, a proper method selection does not guarantee convergence for an arbitrary equation, and the expert system should contain more than methods but proper sets of the parameters. Most of the ODE solvers \cite{hindmarsh2005sundials} are very precise tools yet demanding for the equation form and parameters tuning. 

Solution of partial differential equations is an undeniable traditional topic in the mathematical physics and applications \cite{morton2005numerical}. A wide variety of methods starting from finite-difference schemes through finite element method to modern spectral-like analytical methods are established. In the classical analysis, it is assumed that the operator properties and possible boundary condition types are a priori known since they are defined by the type of the process and the physical nature of the problem that is considered.

The classical finite-difference \cite{thomas2013numerical} and finite-element method \cite{solin2005partial} (FEM) have established area of applicability. For example, FEM is widely used to solve elliptic equations occurring in different areas, for example, mechanics.

Without a doubt, decades of development made FEM the fast method to solve known physics and mechanics-related problems \cite{pavlovic2020geometry}. However, there is no possibility to apply finite-difference and finite element methods to arbitrary equations without significant research of every given problem. Finite-difference methods could be applied to the linear equations. However, it is required to linearize equations, derive a finite difference scheme, and research stability for every PDE solution problem.

Spectral methods \cite{burns2020dedalus} are the most modern analytical and numerical methods for PDE solutions. However, their application to an arbitrary problem is restricted by automatic differentiation on the polynomial decomposition series, restricting the solutions' class.

Arising neural differential operators methods are slightly dependent on a training dataset and require a neural network to be retrained every time a new problem arises \cite{li2020neural}. The recent research shows that combined with the transition to the spectral domain may be promising \cite{li2020fourier} even though it also restricts the applicability to the linear methods if applied to the Fourier specter directly.

To sum up, there is no universal method (and numerous ``No free lunch'' theorems for various numerical algorithms may become evident that it is not possible for differential equations also) that solves every equation fast and precisely. As in the ODE case, a decision support tree could expand every time a new equation and/or boundary condition type appear. Therefore such a system will exchange mathematical physics expert time for the coding and tree forming time. One could name Wolfram Mathematica decision support tree for ODE and PDE as the closest system to the universal. However, it has a drawback of proprietary software and thus is not easily embedded in any algorithm. 

\begin{figure}[ht!]
    \centering
    \includegraphics[width=\linewidth]{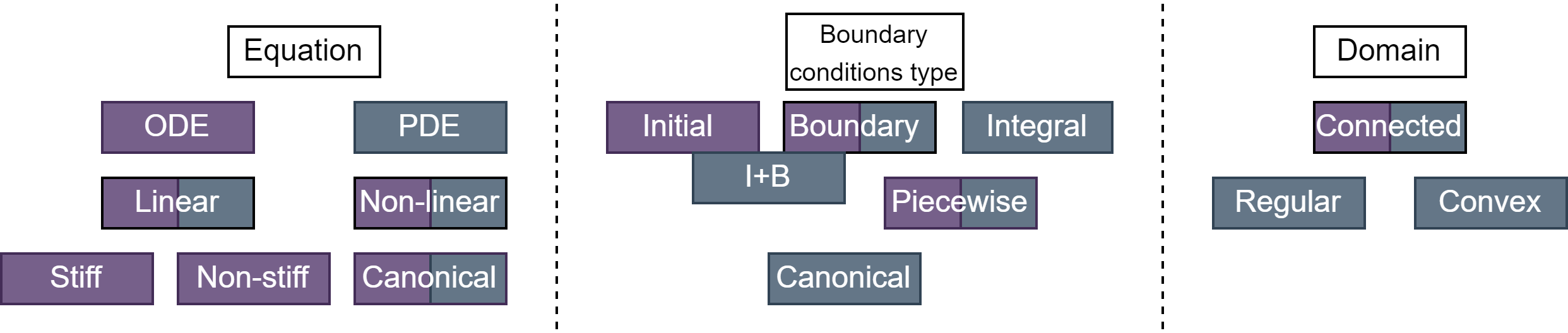}
    \caption{Rough typization of a differential equation solvers.  Typization that is correspondingly may be applied to equation is shown with color (ODE (violet), PDE (blue) or both).}
    \label{fig:solver_scheme}
\end{figure}

The differential equation solution scheme is based on the classification provided in Fig.~\ref{fig:solver_scheme}. We note that real systems contain very extended classification. For every combination  ``equation -- boundary condition type -- domain'' the expert can use a separate solver. For example, the Runge-Kutta method could be chosen in a non-stiff ODE with initial conditions and connected domain case. For stiff equations, more advanced ODE solvers are required. The stiffness assessment, equation and initial condition transformation to the conventional system for Runge-Kutta-like methods are left out of the scope. We emphasize that the equation transformation for the general equation type may be complicated and usually done expertly. In the PDE solution case, we see the same picture. For canonical elastic media equations with arbitrary connected domains, the FEM method family is usually used. Other types of canonical equations usually have finite-difference solutions.

Another critical question is what should be done if the equation cannot be classified using Fig.~\ref{fig:solver_scheme}. We obtain various non-canonical equations with non-canonical boundary conditions during the equation discovery process. Whenever a new equation appears, we should extend classification and add a new solver or apply all available solvers to obtain a possibly incorrect solution. Thus, it is urgent to solve equations in an automated and unified manner to compare solutions with observational data. The less supervised the solution algorithm is, the more viable equation discovery is. Despite being powerful tools in the expert's hands, most of the methods described above are usually applied manually for every given problem. Thus, the topic of an automated arbitrary differential equations solution arises in the literature \cite{rackauckas2020universal}.

The second layer of a problem is the programming language choice. Historically, the solvers are programmed in a performance manner, i.e., using C, C++, and Fortran. Namely, well known solvers ODEPACK (in Fortran) \cite{hindmarsh1983odepack,hindmarsh2019odepack} and solver included in C++ Boost libraries \cite{nakariakov2013boost} are widely used in science. However, modern machine learning and PDE discovery algorithms use Python as a standard, meaning that the machine learning tools usage contradicts the classical approach to differential equation solution. As a compromise, several methods were done using the Julia language \cite{rackauckas2019diffeqflux}. Currently, the Julia language is not so popular among data scientists. In the Python language are available established solvers (as an example, scikit-learn solver and the PyCheb package \cite{liu2016functional}). As an alternative approach, meaning universal solver, possibly single software \cite{lu2021deepxde} could be named. 

In the paper, we describe the automated algorithm of the PDE solution using neural networks and its open-source realization, which reduces computational time and may work with the equation discovery algorithm done in Python. A similar approach is already described in several references and is mainly known as PINN \cite{raissi2019physics}. However, we aim only for the PDE solution part, and thus we omit the physics analogy and replace it with Sobolev space optimization, which is more suitable for general processes.

In contrast to the methods described in the literature, we aim to combine the possibility of solving a wide class of equations and a high level of automatization of the process. It means that we, without the help of the expert, try to combine classical numerical methods and neural networks to obtain the field that approximates the equation's solution without the algorithm parameters changing. The obtained approximation may not have the best quality compared to the classical solvers. The main advantage is that the solution is obtained for every equation despite its type. We use the Sobolev space norm to measure the solution's ``quality''. We validate the algorithm on classical ODEs: Legendre equation and Panleve transcendents and PDEs: wave equation, heat equation, Korteweg-de Vries equation.

The paper is organized as follows:  Sec.~\ref{sec:problem_statement} contains the definitions and algorithm description used in the article, Sec.~\ref{sec:numerical_exps} contains the application of the given algorithm to particular PDEs, Sec.~\ref{sec:conclusions} outlines the paper and proposes the directions for the future work.

\section{Problem statement}
\label{sec:problem_statement}

From the mathematical point of view, the information about the equation during the discovery process is minimal since the algorithms do not use any a priori assumptions about the data governing process. On example of time and single-dimensional space domain, we solve the boundary PDE problem defined on a subdomain $(x,t) \in \Omega \subset 	\mathbb{R} ^2$ with a boundary $\partial \Omega$ in form Eq.~\ref{eq:PDE_problem}. 

\begin{equation}
\begin{array}{cc}
Lu=f\\
bu=g 
\end{array}
\label{eq:PDE_problem}
\end{equation}

In Eq.~\ref{eq:PDE_problem} we assume that the differential operator $L$ and the boundary operator $b$ and the arbitrary functions $f,g$ are defined such that the boundary problem is correct. We do not know any a priori information on the type of $L$ (it may be non-linear, of arbitrary order, with variable coefficients, order and coefficients are subject to change during the discovery process). 

However, we may fix the form of the operator $b$. In classical equation discovery algorithms, boundary conditions are not used, and the equation is found using domain interior. If the dimensionality of data and operator agrees, we may impose classical Dirichlet conditions (function values on the boundary are fixed and taken from data for every equation) to solve a discovered equation. However, the order of operator $L$ is usually not constant in the discovery process, and we may have to take non-conventional conditions, such as the data within the domain interior, to make the problem well-posed. It is worth mentioning that the problem posedness is not used in the algorithm, and theoretically, we may solve under and over-defined problems. In this case, boundary conditions are satisfied in an ``averaged'' manner.

Hypothesis: for such a problem setup, it is more natural to use the optimization methods to obtain the approximation parameters that solve the equation instead of a straightforward solution.

\section{Proposed approach}

\subsection{Theoretical formulation}

Most of the numerical methods assume that the solution field is found in a discrete subset of $\Omega$ in form of the mesh function. It means that for two-dimensional equation we have following function representation:

\begin{equation}
\begin{array}{cc}
\bar{u}=\{u(x^{(i)},t^{(i)}), i=1,2...,n\}\\
\forall i \, (x^{(i)},t^{(i)}) \in \Omega
\end{array}
\label{eq:field_discretization}
\end{equation}

We emphasize that the approach will work for higher dimensions. Without loss of generality, we assume that the field discretization $X\allowbreak=\allowbreak\{x^{(i)}\allowbreak,\allowbreak t^{(i)}\}\subset\Omega$ is fixed during the process of PDE solution. For the experiments, we use a uniform mesh. However, the discretization for the method described below could be chosen arbitrarily.

We formulate a minimization problem to find the solution field as Eq.~\ref{eq:prec_algorithm_formulation}.

\begin{equation}
    \min \limits_{\bar{u}} || L\bar{u}-f||_i+\lambda ||b \bar{u}-g||_j
    \label{eq:prec_algorithm_formulation}
\end{equation}

In Eq.~\ref{eq:prec_algorithm_formulation} arbitrary norms $||\cdot||_i$ and $||\cdot||_j$ may be chosen. Usually  $l_2$ and $l_1$ norms are taken. Operator $L$ is assumed to be the ``precise'' operator that gives the exact value of the derivative for the original function $u(x,t)$ at the mesh points. We note that Eq.~\ref{eq:prec_algorithm_formulation}, $\lambda$ is an arbitrary chosen constant, which does not influence a resulting solution, only convergence speed if the boundary conditions are correctly defined. In this case, there is no doubt that the solution of the optimization problem converges point-wise to the boundary problem solution.

Since the solution is not known, we use numerical differentiation methods to obtain the values of the $L \bar{u}$ for a given solution approximation candidate $\bar{u}$. In practice, differential and boundary operators are also the approximation that has an error, and the minimization algorithm is the numerical algorithm with its own error. Therefore, the final problem that is solved in the article is formulated as Eq.~\ref{eq:approx_algorithm_formulation}.

\begin{equation}
    \min \limits_{\bar{u}} \left[|| \bar{L}\bar{u}-f||_i+\lambda ||\bar{B} \bar{u}-g||_j\right]\Big|_X
    \label{eq:approx_algorithm_formulation}
\end{equation}

In Eq.~\ref{eq:approx_algorithm_formulation} $\bar{L}$ and $\bar{b}$ are the approximate differential and boundary operators (meaning that the derivatives are replaced with the approximations), $X$ is the discretization and $\bar{u}$ are taken accordingly the given grid point. It should be emphasized that the choice of the operator $\bar{L}$  approximation should be considered as a separate problem.

\subsection{Numerical realization}

 Following questions should be answered to solve the problem using any numerical PDE solution method: 
 
 \begin{itemize}
     \item how is the function represented;
     \item how do we take the derivative;
     \item how do we obtain solution approximation parameters.
 \end{itemize}
 
\begin{figure}[h!]
    \centering
    \includegraphics[width=\linewidth]{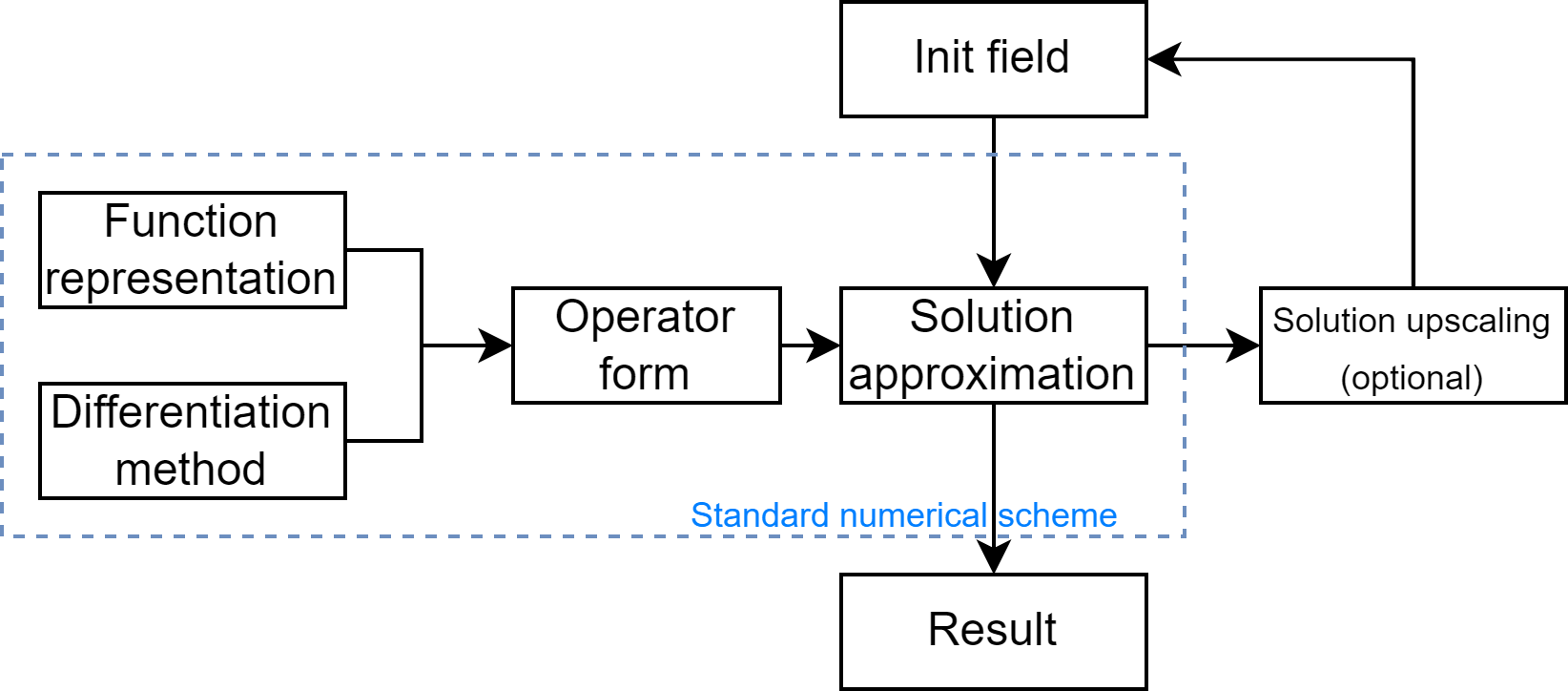}
    \caption{General differential equation solution numerical scheme method}
    \label{fig:PDE_numerical_solution_scheme}
\end{figure}

 The numerical solution scheme is shown in Fig.~\ref{fig:PDE_numerical_solution_scheme}.

As an example for finite-difference scheme we represent the function as the values at the discrete grid points using matrix or more general $n$-dimensional array representation:

\begin{equation}
    \bar{u}=\left( 
    \begin{array}{ccc}
   {{u}_{1,1}} & \cdots  & {{u}_{1,n}}  \\
   \vdots  & \ddots  & \vdots   \\
   {{u}_{m,1}} & \cdots  & {{u}_{m,n}}  \\
    \end{array}
\right)
\label{eq:finite_difference_matrix}
\end{equation}

In this case, numerical differentiation is a condition that binds adjacent nodes and forms the system of linear equations. For example, we can use the first-order finite-difference scheme (approximation order is $O(h)$, where $h$ is the uniform grid step in discretizing the given dimension). We use both forward and backward for boundaries in the form Eq.~\ref{eq:forward_scheme_first}.

\begin{equation}
   \begin{array}{cc}
    u'_f(x) = \frac{u(x+h)-u(x)}{h} \\
    u'_b(x) = \frac{u(x)-u(x-h)}{h}
   \end{array}
   \label{eq:forward_scheme_first}
\end{equation}

We note that in Eq.~\ref{eq:finite_difference_matrix} and in Eq.~\ref{eq:forward_scheme_first} different notations for grid position are used. For the interior points we use scheme Eq.~\ref{eq:center_scheme_first} as more stable and higher-order.

\begin{equation}
u'_c(x) =\frac{1}{2}(u'_f(x) + u'_b(x)) =\frac{u(x+h)-u(x-h)}{2h}
   \label{eq:center_scheme_first}
\end{equation}

As the third component of the numerical PDE solution algorithm, schemes Eq.~\ref{eq:forward_scheme_first}-Eq.~\ref{eq:center_scheme_first} together with the boundary conditions are used to form the system of $n \times m$ equations to find the values at the grid points. We note that such simplicity is not typical for the finite-difference schemes, and usually, the expert is required to form the finite-difference schemes and solutions.

All three components are required to determine most of the numerical PDE solution algorithms. To make the process of PDE solution in a more automated manner, we use machine learning models.  

Meaning that as the first component of the numerical solution we try to approximate solution $u(x,t)$ of an equation $Lu=f$ with continuous parametrized function $\bar{u}(x,t;\Theta) : R^2 \to R$ which is represented by the machine learning model. The parameter set $\Theta=\{\theta_1,...\theta_N \}$ is an arbitrary set that determines the pre-defined function form, as a simplest example $u(x,t; \Theta)=\theta_1 x+ \theta_2 t+ \theta_3$ may be the linear regression. 

As the second component we use finite-difference schemes Eq.~\ref{eq:forward_scheme_first}-Eq.~\ref{eq:center_scheme_first}. We explicitly build the finite difference schemes and combine them to a complete operator for the higher-dimensional derivatives to speed up the computation. However, the finite differences are not used to bind the values in the adjacent grid points. Instead, we find an approximation of whole tangent spaces representing required order derivatives at the given point.

To find the parameter set $\Theta=\{\theta_1,...\theta_N \}$ we use the formulation of the problem Eq.~\ref{eq:approx_algorithm_formulation} in form Eq.\ref{eq:ML_optimization formulation}.

\begin{equation}
    \min \limits_{\Theta} \left[|| \bar{L}\bar{u}(x,t;\Theta)-f||_i+\lambda ||\bar{B} \bar{u}(x,t;\Theta)-g||_j\right]\Big|_X
    \label{eq:ML_optimization formulation}
\end{equation}

In this case, we use the so-called ``mesh-free method". For the neural network, the main difference is that the parameter $h$ may be chosen arbitrarily without connection to the discretization grid $X$. Such solutions are usually referred to as ``mesh-free". 

We note that schemes Eq.~\ref{eq:forward_scheme_first}-Eq.~\ref{eq:center_scheme_first} in the finite-difference analysis is proven to converge in grid points thus taking $h$ same as the grid resolution for schemes is somewhat optimal. We note that there exist different schemes such as $x+1/2 h$ that are not considered for brevity.

To sum up, we encode every operator to use it in the neural network training process. The same procedure is done for the boundary conditions. Starting from the arbitrary neural network's weights set $\Theta_{init}$, we use the optimization algorithm to obtain the optimal set of the parameters $\Theta_{opt}$ that minimizes the difference between the applied operator to the approximated by the neural network function $\bar{L} \bar{u}$ and function $f$ over all discretization points $X$. Additionally, we introduce the difference between the applied boundary operator and function $g$. 

It should be noted that such a learning process differs from the classical neural network training process, where $L^p$ norm convergence is usually considered. In this case, the neural network convergence rather in a Sobolev space \cite{czarnecki2017sobolev}.

\subsection{Modular approach}

\begin{figure}[h!]
 \centering
 \includegraphics[width=\linewidth]{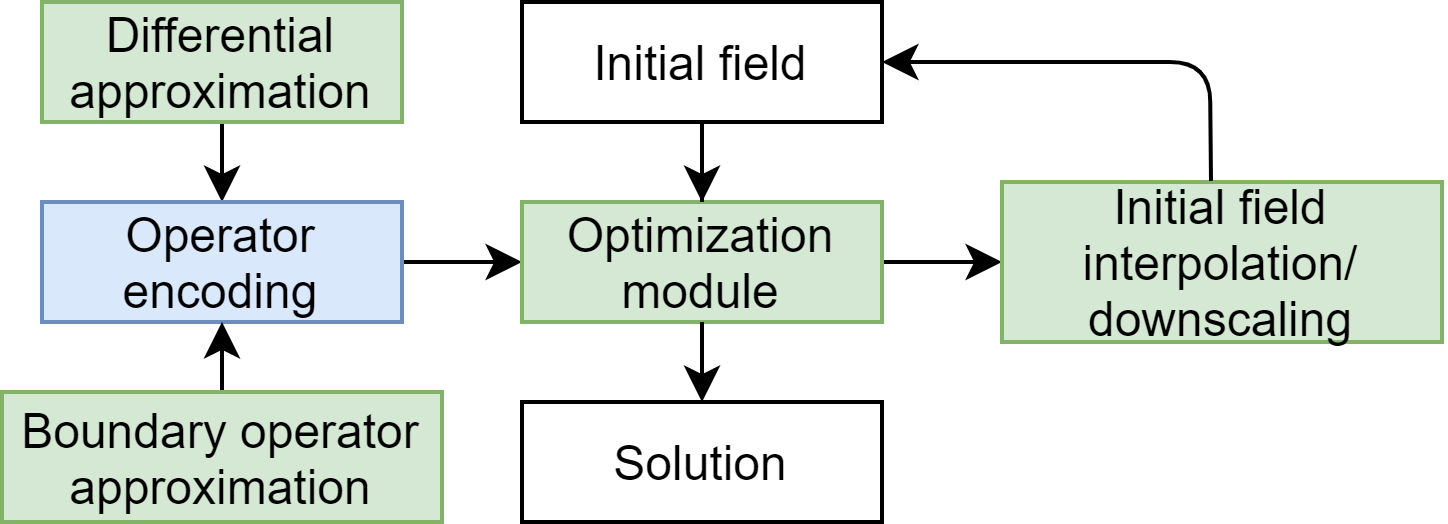}
\caption{Module structure of the solver scheme green parts are replaceable.}
\label{fig:modules}
\end{figure}

To note that the proposed scheme is not unique, we propose the module structure shown in Fig.~\ref{fig:modules} of the resulting solver.

As an alternative, we can use two components of finite difference scheme Eq.~\ref{eq:finite_difference_matrix} and Eq.~\ref{eq:forward_scheme_first}-Eq.~\ref{eq:center_scheme_first} for differential approximation and boundary operator approximation. However, the system is not solved, but matrix values at the grid nodes are optimized such that they satisfy the finite-difference scheme and boundary conditions using the Eq.~\ref{eq:approx_algorithm_formulation}. Such an approach trades quality and continuity to the convergence speed for lower dimension problems.

For neural network training, we use gradient descent. However, the algorithm may be replaced with an evolutionary optimization to make the convergence process more stable due to the possibility of averaging stochastically obtained results. Theoretically, the convergence process may be more stable if proper mutation and cross-over operators are used.

Our experiments show that the problem formulation Eq.~\ref{eq:approx_algorithm_formulation} is more important than the particular realizations of the field approximations, numerical differentiation realization, optimization algorithm, and initial field interpolation or neural network solution upscaling. Thus, we mark corresponding modules in the scheme as replaceable.

\subsection{Caching of approximate models}
\label{sec:cache}

We use initial field interpolation to achieve faster and better convergence. This module replaces the ``initial field interpolation'' module shown in Fig.~\ref{fig:modules}. The effect of the initial guess for the optimization is shown in Sec.~\ref{sec:numerical_exps}. It is realized as a ``cache'' of models for neural networks. As the initial step for every algorithm run, we search in the library the model of the same or another architecture with the lowest Sobolev space norm (sum over all grid points $X$ of the functional in Eq.~\ref{eq:ML_optimization formulation}) for the given equation and boundary conditions. If the model has another input or output dimension, we change the input or output layer and compute the norm. If the models' architectures are different, we train the input architecture on values of a ``cached'' one.

After the algorithm is stopped, the weight of the neural network and the optimizer state (gradient value and related gradient parameters) are saved for further use. Thus, when the same equation is solved again, the shortest possible time is used to obtain the approximation. The use of the caching technique makes the optimization process exploit the existing solution, and thus the optimization tends to converge to the same local optima. The network weights are perturbed every time the best model is taken from the cache to avoid the locality of the solution.

As a result, the proposed approach may solve ODE and PDE similarly without the parameters changing. Even though the parameter's tuning may reduce the optimization time, the overall quality solution remains the same for the broad parameters' range. Such an approach does not challenge the classical methods. On the contrary, the solution in the most challenging cases may be incorrect, but it allows one to compare two equations during the discovery process without stopping the algorithm due to the inappropriate equation for solution errors.

\section{Numerical experiments}
\label{sec:numerical_exps}

The following experiments show the broad range of equations that could be solved with single neural network architecture and algorithm hyperparameters set. Every experiment and picture is supported by repository\footnote{\url{https://github.com/ITMO-NSS-team/torch_DE_solver/tree/main/examples}} with code and experimental data. Experiments show:

\begin{itemize}
    \item that cache allows converging faster;
    \item  that adding points to the grid leads to a better solution;
    \item that the error between the exact and obtained solutions is negligible for equation discovery application.
\end{itemize}

We note that using neural networks does not give the possibility to reproduce singularity points. Therefore all equations are considered in variables' range where no singularities are contained.

As the exact solution, we chose the Wolfram Mathematica 13 \texttt{DSolve} or \texttt{NDSolve} output if it is not stated otherwise.

\subsection{Ordinary differential equations}
\label{sec:exp_ODE}

Ordinary differential equations may be used to determine possible class functions obtained as a solution. This subsection considers two different sets of equations - the Legendre equation and the Panleve transcendents.

\subsubsection{Legendre equation}
\label{sec:exp_Legendre}

Legendre equation may determine which the maximal Taylor series decomposition order may be obtained as a solution. In this section, we consider the problem in form Eq.~\ref{eq:Legendre_eqn}. The solution to the problem is a Legendre polynomial of degree $n$.

\begin{equation}
    (1-t^2)u''(t)-2tu(t)+n(n+1)u=0
    \label{eq:Legendre_eqn}
\end{equation}

Boundary conditions in form Eq.~\ref{eq:Legendre_bcnds} are used.

\begin{equation}
\begin{array}{cc}
    u(0)=L_n(0)\\
    u'(1)=\frac{d L_n(t)}{d t} \Big|_{t=1}
\end{array}
    \label{eq:Legendre_bcnds}
\end{equation}

In Eq.~\ref{eq:Legendre_bcnds} $L_n(t)$ is a Legendre polynomial of degree $n$.

\begin{figure}[ht!]
    \centering
    \includegraphics[width=\linewidth]{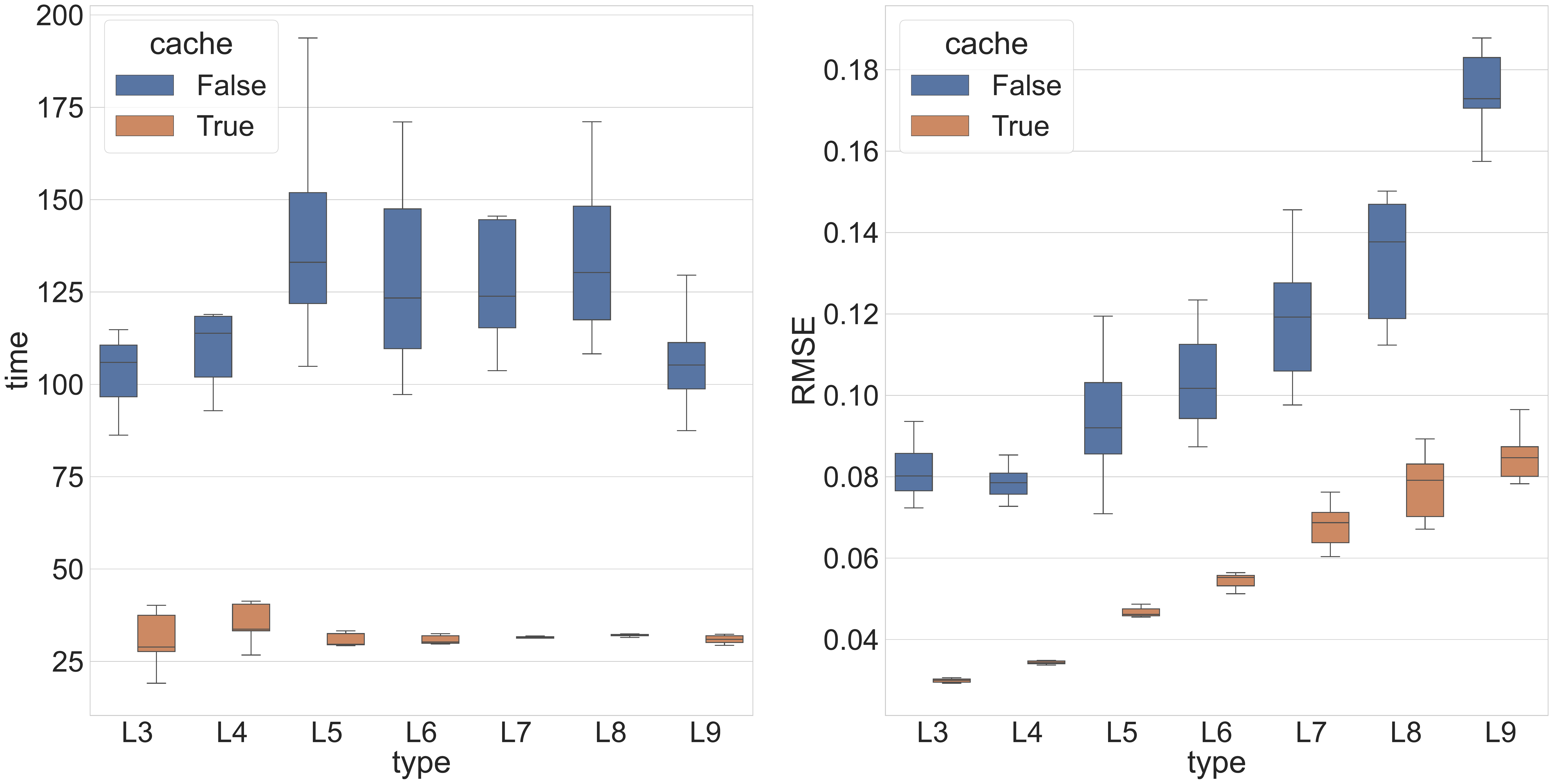}
    \caption{Legendre polynomial solution time  in seconds (left) and solution error (right). Type is the order of polynomial $n$ in Eq.~\ref{eq:Legendre_eqn}, with different colors  (\texttt{cache=false}) and (\texttt{cache=true}) cases are shown}
    \label{fig:legendre_exps}
\end{figure}

First set of the experiments is the learning of the neural network parameters for $n=3,...,9$ for $100$ uniformly taken points from a range $t\in [0;1]$. The time of optimization was recorded without (\texttt{cache=false}) and with (\texttt{cache=true}) caching technique (see Sec.~\ref{sec:cache}) for $10$ runs are shown in Fig.~\ref{fig:legendre_exps} (left).

As seen from Fig.~\ref{fig:legendre_exps} (left), the reasonable initial guess, as expected, makes the optimization converge faster. As a drawback, such an approach makes the optimization more ``stiff'' and thus, in some cases, it may get stuck in the local minima. The initial guess parameters are perturbed with a noise of a small magnitude to mitigate stiffness partially.

The root mean square errors (RMSE) for the same setup are shown in Fig.~\ref{fig:legendre_exps} (right). The error is computed using the analytical solution - Legendre polynomial values of corresponding order in a range $t \in [0,1]$ (the half range is taken due to symmetry property) for every grid point. Since the maximal value of Legendre polynomial on range $t \in [0,1]$ is 1, RMSE could be interpreted as the absolute error.

In Fig.~\ref{fig:legendre_exps} (right) we see that the locality of the solution described in Sec.~\ref{sec:cache} appears when the cache is used. Namely, the error dispersion is lower when an initial guess is used. It allows the algorithm to find the solution with a lower Sobolev space norm as a positive effect. Therefore, the initial guess allows faster convergence, and the solution is likely to have a better norm.

Overall, the ability of the solver to converge towards a Legendre polynomial solution means that the solver can converge towards any analytical solution (a solution that may be represented in the form of the Taylor series). We may also answer on a maximal Taylor series decomposition order. We obtain at least the ninth term in decomposition with a good (less than 10 \% error) precision for $100$ points without parameter change and optimization time increase.

\begin{figure}[ht!]
    \centering
    \includegraphics[width=\linewidth]{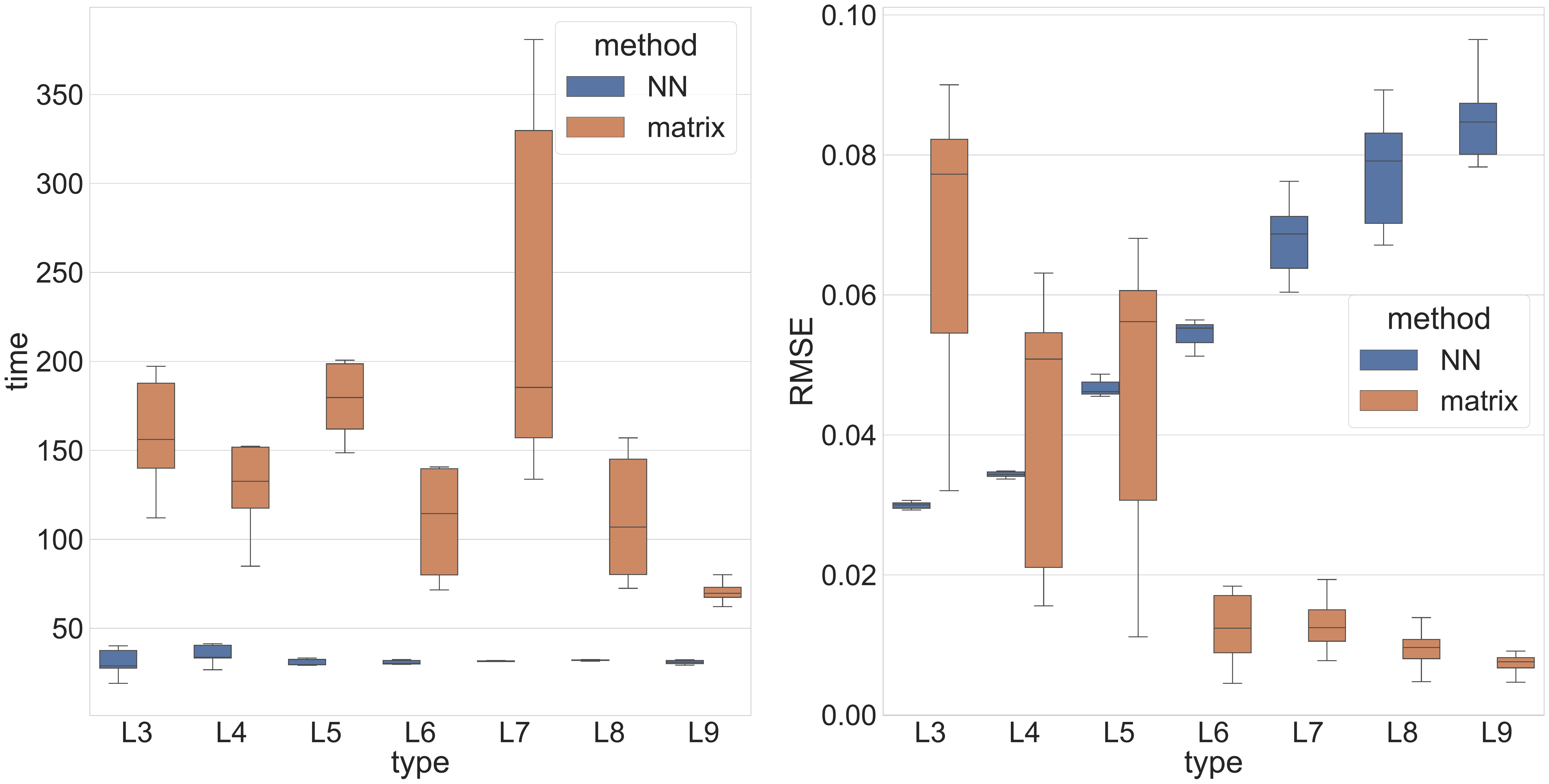}
    \caption{Legendre polynomial solution time (left) and solution error (right). Type is the order of polynomial $n$ in Eq.~\ref{eq:Legendre_eqn}, different colors - different approximations - neural network-based and matrix-based}
    \label{fig:Legendre_error_time_mat}
\end{figure}

To reduce error and show different module from Fig.~\ref{fig:modules}. We show the solution using the matrix differentiation prototype. The results for same grid setup for Legendre equation Eq.~\ref{eq:Legendre_eqn} are shown in Fig.~\ref{fig:Legendre_error_time_mat}.

From comparison matrix and neural network approximation, we see that whereas matrix has a lower error, it has higher computation time. Thus, additional techniques such as consequent field up-sampling or another result ``caching'' are required. The matrix-based algorithm is, however, out of the paper's scope. Therefore further experiments are conducted using neural network approximation only.

\subsubsection{Panleve transcendents}
\label{sec:exp_Panleve}

\begin{figure}[ht!]
    \centering
    \includegraphics[width=\linewidth]{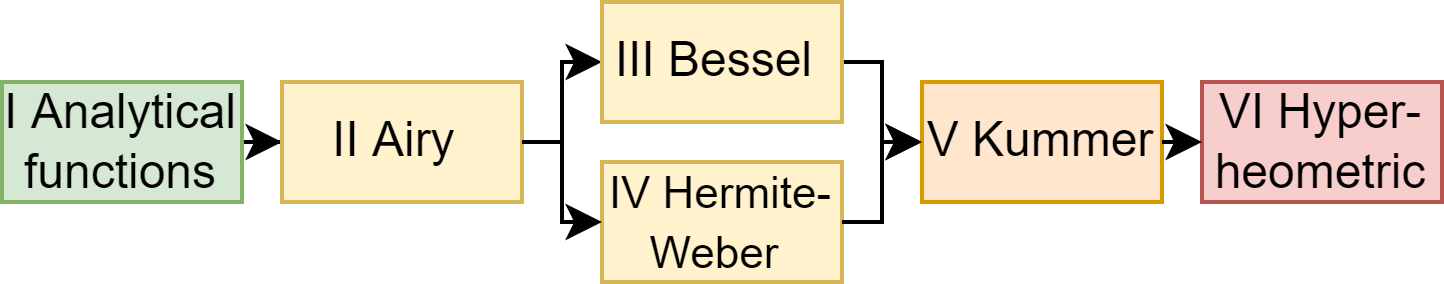}
    \caption{Scheme of solution of Panleve transcendent ``complexity''. Next class of functions contains previous. The color is predicted relative ``performance'' of a conventional  solver for every transcendent.}
    \label{fig:Panleve_scheme}
\end{figure}

Panleve transcendents are a series of differential equations. Each has a different special function class in the solution. Meaning that the higher order of transcendent moves a step closer to the general hypergeometric function. This subsection considers different aspects of the Panleve transcendents solution as an essentially non-linear ODE with variable coefficients. The scheme of functions that form the general solution of a given Panleve transcendent is shown in Fig.~\ref{fig:Panleve_scheme}.

Exact initial and boundary problems are placed in \ref{app:Panleve} since the form of the equation is not important for further experiments. We note that the value range is taken such that the solution does not contain any singularity points.

Whereas the Legendre equation appears as a relatively simple linear equation with variable coefficients, the Panleve transcendents are significantly non-linear and have wider solution space than the polynomial. Additionally, the maximal sequential number of transcendent allows showing which class of function solver can reproduce. A hypergeometric function is the broadest possible class for a real-valued ODE solution.

We note that the solver parameters and the neural network model are the same as for the Legendre polynomial in Sec.~\ref{sec:exp_Legendre} in most Panleve transcendents experiments (with possibly variable stop-criterion, meaning that in some cases optimization is stopped earlier or later depending on the equation).

\begin{figure}
    \centering
    \includegraphics[width=\linewidth]{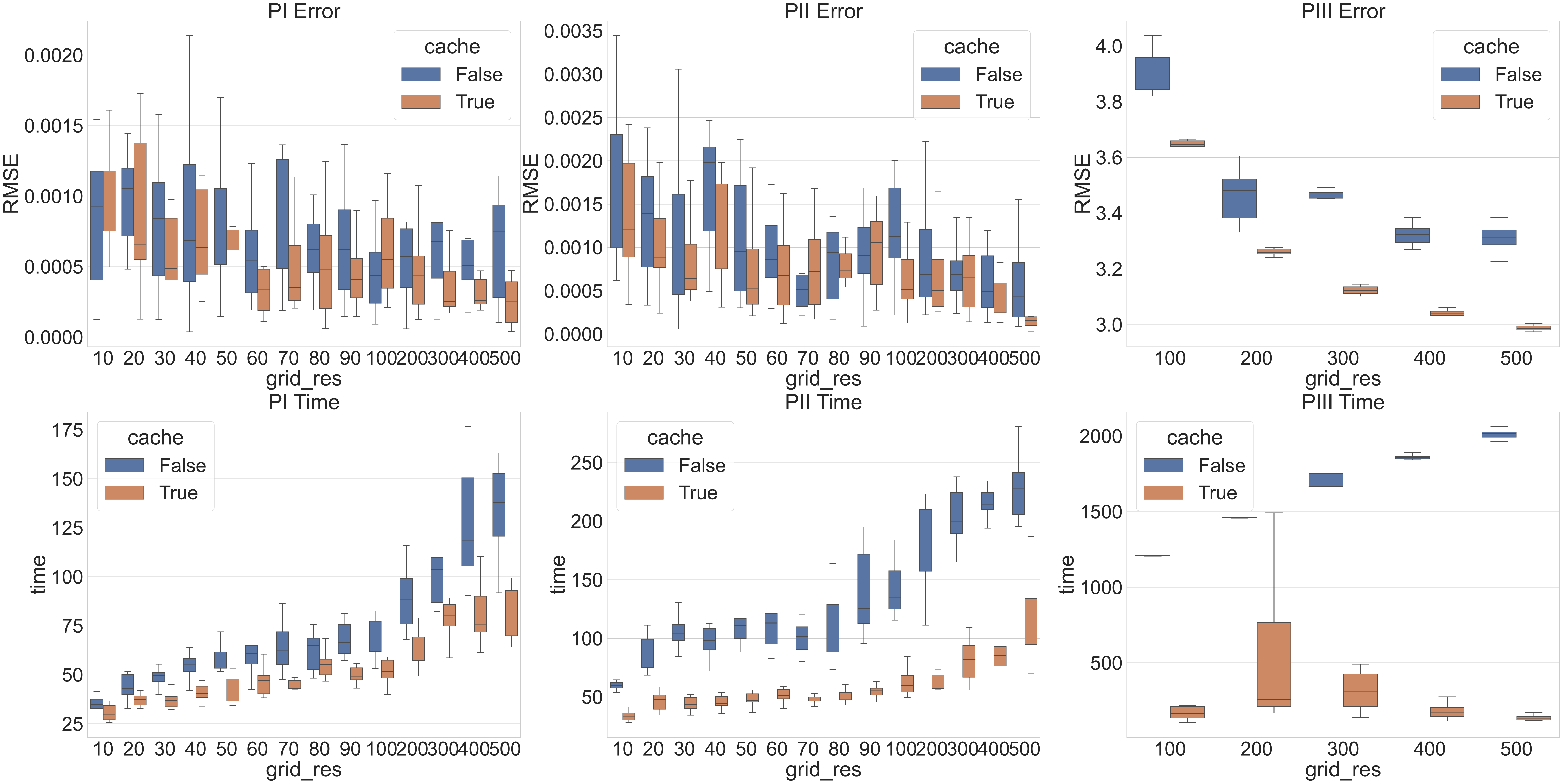}
    \caption{A summary of Panleve I-III expreriments runs, different colors are the runs with and without initial guess}
    \label{fig:Panleve_exps_I_III}
\end{figure}

For the first three Panleve transcendents, we repeat the same experimental setup as in Sec.~\ref{sec:exp_Legendre}. However, this experiment series is used to show convergence by changing the amount of uniformly taken points $grid\_res$ from a value range $t\in T$. Error and optimization time distribution are shown in Fig.~\ref{fig:Panleve_exps_I_III}. Error is computed using the Wolfram Mathematica 13 numerical solution on optimization grid points.

As seen from Fig.~\ref{fig:Panleve_exps_I_III} using proper initial weights distribution reduces optimization time and possibly leads to a lower error. 

\begin{figure}[h!]
    \centering
    \includegraphics[width=\linewidth]{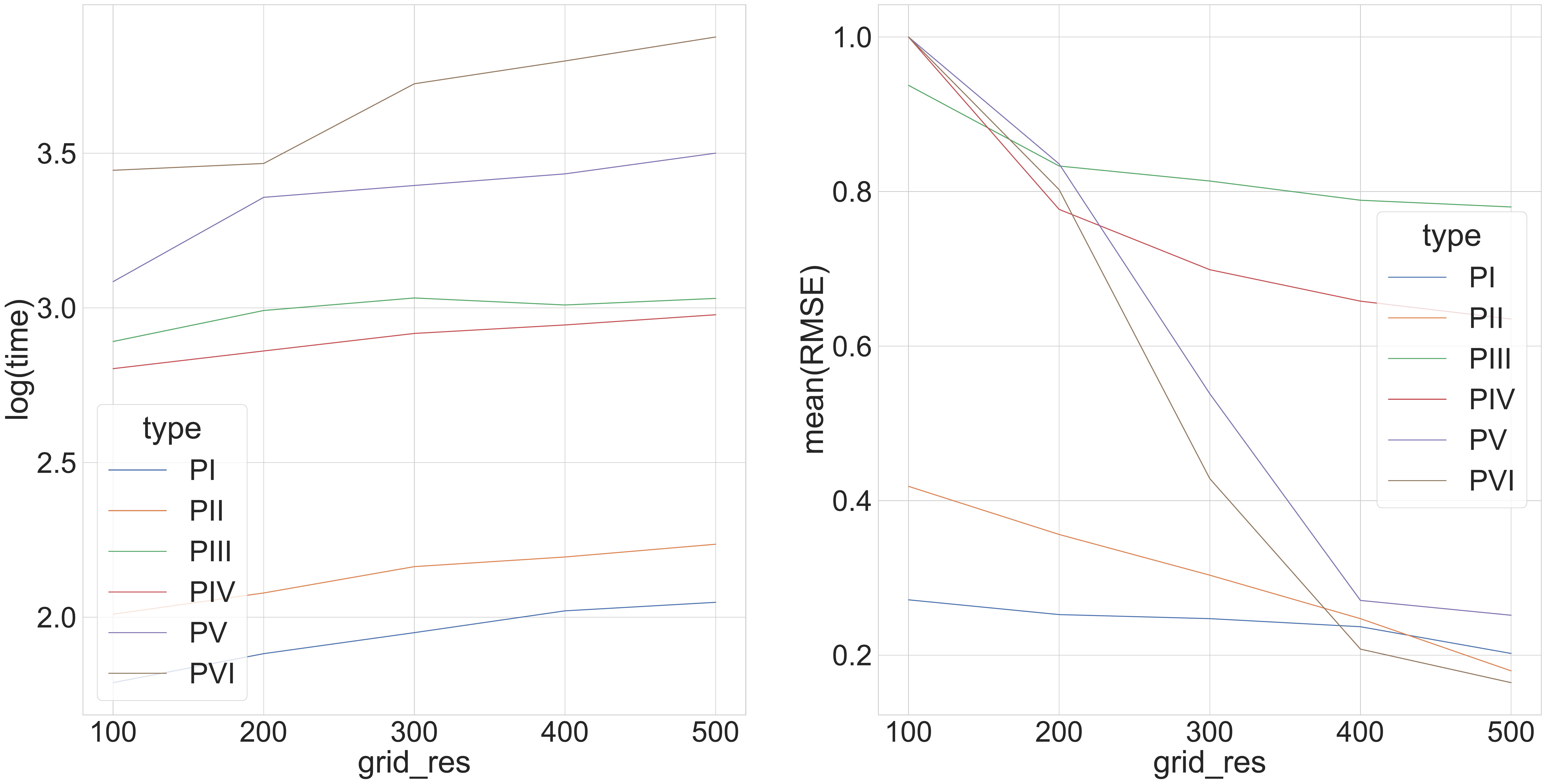}
    \caption{Mean time (in seconds) in logarithmic scale of optimization for each Panleve transcendent (left) and mean error ratio (right) (1.0 is a maximal error for each equation), maximal error for PI-PII is at the $10$ grid points, which is not shown}
    \label{fig:Panleve_exps_mean}
\end{figure}

More complex transcendents have increased solving time as shown in Fig.~\ref{fig:Panleve_exps_mean} (left) (we emphasize that the logarithmic scale for time is used). 

We note that the optimization time mostly depends on the number of terms of the equation. To assess the influence of the solution complexity, we should use the set of the equation of a similar number of terms with different solution classes, which is nearly impossible.

For all six Panleve transcendents mean (for PIV-PVI only one run was performed) error over all experiments using cache is shown in Fig.~\ref{fig:Panleve_exps_mean} (right). We emphasize that error is normalized on the maximum error achieved during all experiments for all grid points.

In summary, algorithm error ``converges'' towards a solution for every transcendent when there are no discontinuity points. Thus, we may conclude that the optimization problem always has a solution close to the true PDE solution in the range where the solution is analytic. However, such a statement requires more rigorous proof out of the paper scope.

\subsection{Partial differential equations}

\begin{figure}[ht!]
    \centering
    \includegraphics[width=\linewidth]{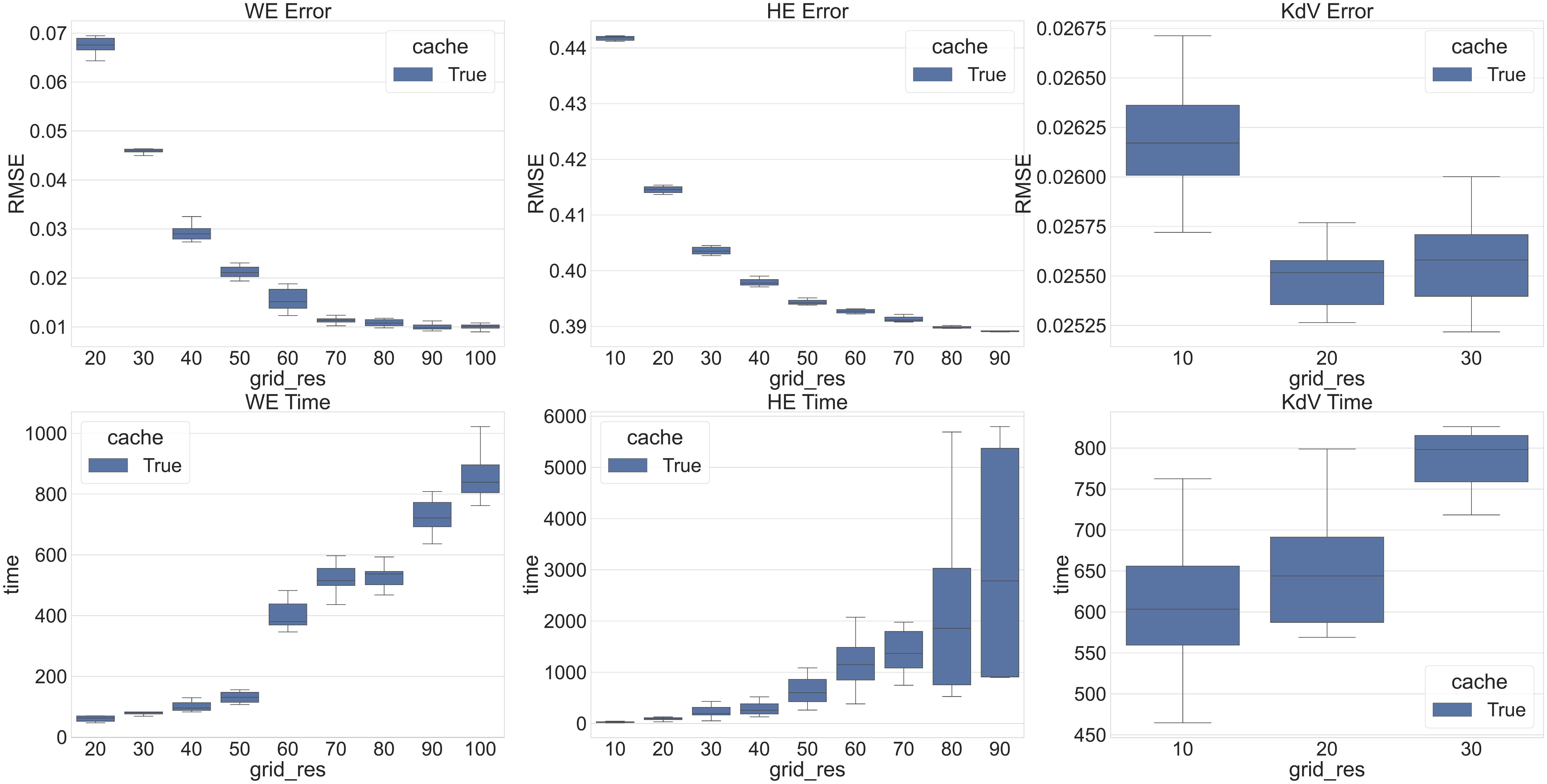}
    \caption{Summary of all PDE experiments: upper row - errors, lower - time. From left to right - wave equation, heat equation, Korteweg-de Vries equation \label{fig:PDE_exps}}
\end{figure}

Unlike ODEs, PDEs contain an interconnection between differentials with respect to the different variables. Therefore, the ability of the ODE solution does not imply the ability to solve PDE. In this section, two canonical equations - wave and heat equations- are considered a ``lower'' bound of the complexity, and the Korteweg-de Vries equation is considered a more complex example of a non-linear equation. We note that it is impossible to pick an ``upper'' bound of PDE solution due to the non-classicability of the characteristic surface of the equation starting from order three.

\subsubsection{Wave equation}
\label{sec:exp_wave}

We try to assess the convergence of the algorithm the solution of the wave equation with boundary conditions in form Eq.~\ref{eq:WE_formulation}

\begin{equation}
    \begin{array}{cc}
         \frac{\partial^2 u(x,t)}{\partial t^2}-\frac{1}{4} \frac{\partial^2 u(x,t)}{\partial x^2}=0 \\
         u(0,t)=u(1,t)=0 \\
         u(x,0)=u(x,1)=\sin(\pi x) \\
         (x,t) \in [0,1] \times [0,1]=\Omega
    \end{array}
    \label{eq:WE_formulation}
\end{equation}

We use the formulation Eq.\ref{eq:prec_algorithm_formulation} to obtain the solution of the equation for ten runs for consequently increasing the number of points in discretization from $10 \times 10$ points in $\Omega$ (since the mesh is assumed uniform, it is equal to $h=1/(N-1)=1/9$, where $N$ is the number of points for time and space dimensions, i.e., we take 10 points in the range $[0,1]$ including boundaries) to $100 \times 100$ points with the step of 100 points. As a result of ODE experiments, we show only the ``cache'' version, meaning that every time we start with the best possible initial neural network weights (see Sec.~\ref{sec:cache}).

We take an analytical solution from the Wolfram Mathematica 13 software as the exact solution. The solution has the analytical form and is taken at the grid points for every grid used in the optimization process. We record the optimization time, and the root mean square error between (RMSE) Wolfram Mathematica solution and the proposed algorithm solution on the same grid. The time and error distribution for 10 runs are shown in Fig.~\ref{fig:PDE_exps} (left).



\subsubsection{Heat equation}
\label{sec:exp_heat}

The second equation is a typical equation of a parabolic type - heat equation. To demonstrate that the algorithm converges in case of an incorrectly posed boundary problem, we use the following formulation shown in Eq.~\ref{eq:HE_formulation}

\begin{equation}
    \begin{array}{cc}
         \frac{\partial u(x,t)}{\partial t}-\frac{\partial^2 u(x,t)}{\partial x^2}=0 \\
         u(0,t)=500 \\
         \frac{\partial u(x,t)}{\partial x} \Big |_{x=1}=1 \\
         u(x,0)=0 \\
         (x,t) \in [0,1] \times [0,1]=\Omega
    \end{array}
    \label{eq:HE_formulation}
\end{equation}

The boundary problem Eq.~\ref{eq:HE_formulation} for heat equation has an analytical solution in form Eq.~\ref{eq:HE_sln}.

\begin{multline}
        u(x,t)=500+x+\\+\frac{8}{\pi^2} \sum_{k=1}^{k=+\infty} \Big[\exp{(-\frac{1}{4} \pi ^2 t (2 k-1)^2)}  \\ \frac{ \left(250
   \pi  (1-2 k)+(-1)^{k}\right) \sin \left(\frac{\pi  x }{2}
   (2 k-1)\right)}{(2 k-1)^2} \Big]
    \label{eq:HE_sln}
\end{multline}

We use the same grid setup for the experiments as for the wave equation in Sec.~\ref{sec:exp_wave}. Namely, ten runs from $10 \times 10$ points to $100 \times 100$ uniformly taken from $ \Omega =[0,1] \times [0,1]$. The error in this case is computed using the analytical solution Eq.~\ref{eq:HE_sln} with 100 first terms in sum taken. The error and time plots are shown in Fig.~\ref{fig:PDE_exps} (middle).



\subsubsection{Korteweg-de Vries equation}
\label{sec:exp_KdV}

 To show a more sophisticated PDE solution, the Korteweg-de Vries equation (Eq.~\ref{eq:KdV}) was used. 

\begin{equation}
\label{eq:KdV}
u_t + 6 u u_x + u_{xxx} = f(x,t)
\end{equation}

Following forcing, initial and boundary conditions were applied as shown in Eq.~\ref{eq:KdV_bc}.

\begin{equation}
\label{eq:KdV_bc}
\begin{array}{cc}
f(x,t)=\cos t \sin x\\
u(x,0)=0\\
\left[u_{xx}+2 u_x+u \right] \Big|_{x=0}=0\\
\left[2 u_{xx}+u_x+3 u \right] \Big|_{x=1}=0\\
\left[5 u_x+5 u \right] \Big|_{x=1}=0\\
\end{array}
\end{equation}

Due to the extended computation time, the KdV equation solution was tested on a $10 \times 10$, $20 \times 20$ and $30 \times 30$ uniformly taken points from a range $x \times t \in [0,1]\times[0,1]$. The results of 10 consequent runs for every experiment are shown in Fig.~\ref{fig:PDE_exps} (right). We note that initial solution was obtained within 2000 seconds time range.



Even though the computation time is higher than in the ODE case, the PDE part allows obtaining at least an approximated coarse solution, which could be used in the equation discovery algorithm within a reasonable time. Moreover, the error is below 10 \% of maximal field value at the initial $10 \times 10$ grid.

\section{Conclusions}
\label{sec:conclusions}

The paper proposes a unified numerical differential equation solver based on optimization methods. It has the following advantages:

\begin{itemize}
    \item It can solve PDEs without the involvement of an expert, which is most useful for data-driven equation discovery methods;
    \item It has good precision every for equation discovery application;
    \item It can be easily parallelized;
    \item It has a flexible modular structure. The modules could be replaced to achieve better speed or better precision.
\end{itemize}

 It is also seen that the optimization time is the main drawback of the method's experimental realization. We propose the following optimization speed-up directions:

\begin{itemize}
    \item usage the power of GPU to make the optimization using fast memory and built-in matrix instructions
    \item Better usage of initial approximation
    \item More intelligent use of the numerical differentiation
    \item Usage of joint neural network and matrix-based optimization methods
\end{itemize}

All experimental data and script that allows reproducing experiments are available at the GitHub repository \footnote{\url{https://github.com/ITMO-NSS-team/torch_DE_solver}}.

\section*{Acknowledgments}

This research is financially supported by The Russian Scientific Foundation, Agreement \#21-71-00128.

\appendix

\section{Panleve boundary problems}
\label{app:Panleve}

For all equations $\alpha=\beta=\gamma=\delta=1$.

Panleve I:

\begin{equation}
   \begin{array}{cc}
         u''(t)=6 u(t)^2+t\\
        u(0)=0\\
        u'(0)=0\\
        t \in [0;1]
   \end{array}
   \label{eq:PI}
\end{equation}

Panleve II:

\begin{equation}
   \begin{array}{cc}
        u''(t)=\alpha +2 u(t)^3+t u(t)\\
        u(0)=0\\
        u'(0)=0\\
        t \in [0;1]
   \end{array}
   \label{eq:PII}
\end{equation}

Panleve III:

\begin{equation}
   \begin{array}{cc}
            t u(t) u''(t)=\delta  t-u(t) u'(t)+t u'(t)^2+\alpha  u(t)^3+\\
            +\beta  u(t)+\gamma  t u(t)^4\\
        u(1)=0\\
        u'(1)=0\\
        t \in [0.25;2.1]
   \end{array}
   \label{eq:PIII}
\end{equation}

Panleve IV:

\begin{equation}
   \begin{array}{cc}
    u(t) u''(t)=\beta +2 \left(t^2-\alpha \right) u(t)^2+\\
    +\frac{1}{2} u'(t)^2+\frac{3 u(t)^4}{2}+4 t u(t)^3\\
        u(1)=0\\
        u'(1)=0\\
        t \in [1/4;7/4]
   \end{array}
   \label{eq:PIV}
\end{equation}

Panleve V:

\begin{equation}
   \begin{array}{cc}
         2 t^2 (1-u(t)) u(t) u''(t)=2 \beta+\\
         +2 u(t)^2 \left(\alpha +3 \beta -\delta  t^2+\gamma  t+t u'(t)\right)-\\
         -u(t) \left(6 \beta +3 t^2 u'(t)^2+2 t u'(t)\right)+\\
         +t^2 u'(t)^2-2 u(t)^3 (3 \alpha +\beta +t (\gamma +\delta  t))-\\
         -2 \alpha  u(t)^5+6 \alpha  u(t)^4\\
        u(0.9)=3\\
        u(1.2)=4\\
        t \in [0.9;1.2]
   \end{array}
   \label{eq:PV}
\end{equation}

Panleve VI:

\begin{equation}
\begin{array}{lr}
c_{6,0}+c_{6,1} u(t)+c_{6,2} u(t)^2+c_{6,3} u(t)^3+c_{6,4} u(t)^4+\\
+c_{6,5} u(t)^5-\alpha  u(t)^6+c_{6,6} u(t)u'(t)+c_{6,7} u(t)^2 u'(t)+\\
+c_{6,8} u(t)^3 u'(t)++c_{6,9} u'(t)^2+c_{6,10} u(t) u'(t)^2+\\
   +c_{6,11} u(t)^2 u'(t)^2+c_{6,12} u(t) u''(t)+c_{6,13} u(t)^2 u''(t)+\\
   +c_{6,14} u(t)^3 u''(t)=0 \\
c_{6,0}= -t^3 \beta \\
c_{6,1}=2 \beta  t^2 (t+1)\\
c_{6,2}=-t (\beta -\delta +t (\alpha +\delta +\beta  (t+4)+\gamma  (t-1)))\\
c_{6,3}=2 t (\alpha  (t+1)+\beta  (t+1)+(t-1) (\gamma +\delta ))\\
c_{6,4}=-\alpha +\gamma -\alpha  t (t+4)-t (\beta +\gamma +\delta  (t-1))\\
c_{6,5}=2 \alpha (t+1)\\
   c_{6,6}=(t-1) t^3\\
   c_{6,7}=-t \left(t \left(t^2+t-3\right)+1\right)\\
   c_{6,8}=(t-1) t (2 t-1)\\
   c_{6,9}=-\frac{1}{2} (t-1)^2 t^3\\
   c_{6,10}=(t-1)^2 t^2 (t+1)\\
   c_{6,10}=-\frac{3}{2} (t-1)^2 t^2\\
   c_{6,11}=(t-1)^2 t^3\\
   c_{6,12}=-(t-1)^2 t^2 (t+1)\\
   c_{6,13}=(t-1)^2 t^2 \\
    u(1.2)=u(1.4)=2\\
    t \in [1.2;1.4] 
\end{array}
   \label{eq:PVI}
\end{equation}

%
%
%
 \bibliographystyle{splncs04}
\bibliography{references}

\end{document}